\documentclass[12pt]{article}

\usepackage{amsmath,amssymb,url}
\usepackage{tikz,color,pgf}
\usepackage{longtable} 

\newtheorem{theorem}{Theorem}[section]

\newtheorem{lemma}[theorem]{Lemma}
\newtheorem{proposition}[theorem]{Proposition}
\newtheorem{conjecture}[theorem]{Conjecture}

\newtheorem{example}[theorem]{Example}

\def\cp{\,\square\,}
\newcommand{\proof}{\noindent{\bf Proof.\ }}
\newcommand{\qed}{\hfill $\square$ \bigskip}
\newcommand{\dstart}{\gamma_g}
\newcommand{\sstart}{\gamma_g^\prime}
\newcommand{\w}{{\rm w}}

\title{On Rall's $1/2$-conjecture on the domination game}

\author{Csilla Bujt\'as$\/^{a}$, Vesna Ir\v{s}i\v{c}$\/^{b}$, Sandi Klav\v{z}ar$^{a,b,c}$, Kexiang Xu$\/^{d}$ \\\\
$^{a}$ \small Faculty of Mathematics and Physics, University of Ljubljana, Slovenia\\
$^{b}$ \small Institute of Mathematics, Physics and Mechanics, Ljubljana, Slovenia \\
$^{c}$ \small Faculty of Natural Sciences and Mathematics, University of Maribor, Slovenia\\
$^{d}$ \small  College of Science, Nanjing University of
 Aeronautics \& Astronautics,\\
 \medskip
  \small Nanjing, Jiangsu 210016, PR China\\
 \small {\tt csilla.bujtas@fmf.uni-lj.si},\quad  \small {\tt vesna.irsic@fmf.uni-lj.si}\\
  \small{\tt sandi.klavzar@fmf.uni-lj.si},\quad   \small {\tt kexxu1221@126.com}
}
\date{\today}

\renewcommand{\gg}{\gamma_{\rm g}}

\begin{document}

\maketitle

\begin{abstract}
The $1/2$-conjecture on the domination game asserts that if $G$ is a traceable graph, then the game domination number $\gamma_g(G)$ of $G$ is at most $\left\lceil \frac{n(G)}{2} \right\rceil$. A traceable graph is a $1/2$-graph if $\gamma_g(G) = \left\lceil \frac{n(G)}{2} \right\rceil$ holds. It is proved that the so-called hatted cycles are $1/2$-graphs and that unicyclic graphs fulfill the $1/2$-conjecture. Several additional families of graphs that support the conjecture are determined and computer experiments related to the conjecture described.
\end{abstract}


\medskip\noindent
\textbf{AMS Math.\ Subj.\ Class.\ (2010)}: 05C57, 05C69, 91A43

\section{Introduction}

The domination game on a graph $G$ is played by two players referred to as Dominator and Staller. If Dominator (resp.\ Staller) is the one to start the game, we speak of the D-game (resp.\ S-game). The players alternately select vertices such that at each move at least one vertex is dominated that has not yet been dominated by the set of previously selected vertices. As soon as this is not possible, the game is over; at that point the selected vertices form a dominating set of $G$. Dominator's goal is to reach the end of the game as soon as possible, while Staller has the opposite goal. Assuming that both players play optimally, the number of moves played in the D-game (resp.\ S-game) is a graph invariant  denoted by $\gamma_g(G)$ and named {\em game domination number}  (resp.\ {\em Staller-start game domination number} $\gamma_g'(G)$) of $G$.

The domination game was introduced in~\cite{bresar-2010}. Early influential references on this game include~\cite{bresar-2013, bujtas-2015a, bujtas-2015b, dorbec-2015, kinnersley-2013, kosmrlj-2014}, while from the very extensive recent development on the domination game and its variants we select the papers~\cite{borowiecki-2019, bresar-2019, bujtas-2019, irsic-2019, jiang-2019, ruksasakchai-2019, xu-2018}.  In this paper we are interested in the following conjecture that has been proposed several years ago by D.~Rall, the first published source of it is~\cite[Conjecture 1.1]{james-2019}. Recall that a {\em traceable graph} is a graph that contains a Hamiltonian path.

\begin{conjecture}
\label{conj:Rall}
If $G$ is a traceable graph, then $\gamma_g(G) \le \left\lceil \frac{n(G)}{2} \right\rceil$.
\end{conjecture}

We say that a graph $G$ is a {\em $1/2$-graph} if $G$ is traceable and $\gamma_g(G) = \left\lceil \frac{n(G)}{2} \right\rceil$.  In other words, $1/2$-graphs are the traceable graphs that attain the equality in Conjecture~\ref{conj:Rall}. The existence of families of $1/2$-graphs (see Section~\ref{sec:1/2-graphs}) implies that, if Conjecture~\ref{conj:Rall} holds true, then the asserted bound is best possible.

The paper is organized as follows. In the next section we collect notation and earlier results needed in this paper. In Section~\ref{sec:1/2-graphs} we first recall families of $1/2$-graphs that are already known, and prove that the so-called hatted cycles are $1/2$-graphs as well. In Section~\ref{sec:unicyclic} we prove that unicyclic graphs fulfill Conjecture~\ref{conj:Rall}.  Then, in Section~\ref{sec:additional} we determine several additional families of graphs that support the conjecture. In the subsequent section we report on our computer experiments related to the conjecture.

\section{Preliminaries}
\label{sec:preliminaries}

For a positive integer $k$ we will use the notations $[k] = \{1,\ldots, k\}$ and $[k]_0 = \{0,1,\ldots, k-1\}$. The order of a graph $G$ will be denoted by $n(G)$. A vertex of a graph \emph{dominates} itself and its neighbors; a \emph{dominating set} in a graph $G$ is a set of vertices of $G$ that dominates all vertices in the graph. The cardinality of a smallest dominating set of $G$ is the {\em domination number} $\gamma(G)$ of $G$. If $G$ is a graph and $S\subseteq V(G)$, then a {\em partially dominated graph} $G|S$ is a graph together with a declaration that the vertices from $S$ are already dominated. We will need the following fundamental result on the domination game.

\begin{lemma}{\rm \cite[Lemma~2.1]{kinnersley-2013}}
  \label{lem:continuation}
  {\rm (Continuation Principle)}
  Let $G$ be a graph, and let $A,B\subseteq V(G)$. If $B\subseteq A$, then
  $\dstart(G|A)\le \dstart(G|B)$ and $\sstart(G|A)\le \sstart(G|B)$.
\end{lemma}

Let $P_n' = P_{n+1}|u$ and $P_n'' = P_{n+2}|\{u,v\}$, where $u$ and $v$ are the end-vertices of the path in question. We will also need the following result.

\begin{lemma}{\rm \cite[Lemmas~2.1 and 2.3]{kosmrlj-2017}}
  \label{lem:P'-P''}
  If $n\ge 0$, then
  \begin{eqnarray*}
  \dstart(P'_n) = \dstart(P''_n)&\!\!\!=\!\!\!&\left\{\begin{array}{ll}\lceil \frac{n}{2}\rceil-1;& n \equiv 3\ (\bmod\ 4), \\ \lceil \frac{n}{2}\rceil;& {\rm otherwise.} \end{array}\right. \\
  \sstart(P'_n) = \sstart(P''_n)&\!\!\! = \!\!\! &\left\{\begin{array}{ll}\lceil \frac{n}{2}\rceil+1;& n \equiv 2\ (\bmod\ 4), \\\lceil \frac{n}{2}\rceil;& {\rm otherwise.} \end{array}\right.
  \end{eqnarray*}
\end{lemma}

Define the {\em weighting function} $\w$ of partially dominated paths $P'_{4q+r}$ and $P''_{4q+r}$ with
\[
\w(P'_{4q+r})= \w(P''_{4q+r})= 2q +
  \begin{cases}
    0; &  r=0,\\
    1; & r=1,\\
    \frac{3}{2}; & r=2, \\
    \frac{7}{4}; & r=3.
  \end{cases}
\]

Here is another result that will be applied.

\begin{lemma}{\rm \cite[Lemma 3.2]{dorbec-2019}}
\label{lem:union}
 {\rm (Union Lemma)}
If $F_1, \ldots, F_k$ are vertex-disjoint paths where $F_i = P_{n_i}'$ or $F_i = P_{n_i}''$ for $i \in [k]$ and $n_i \ge 1$, then
\[
\gamma_g'\left(\bigcup_{i = 1}^k F_i \right) \le \left\lceil \sum_{i = 1}^k \w(F_i) \right\rceil.
\]
\end{lemma}

Note that Lemma~\ref{lem:union} clearly remains true under the weaker condition $n_i \ge 0$, $i \in [k]$.

\section{Families of $1/2$-graphs}
\label{sec:1/2-graphs}

In this section we first state which paths and cycles are $1/2$-graphs, then recall that broken ladders are $1/2$-graphs, and end the section by proving that the so-called hatted cycles are also $1/2$-graphs.

\subsection*{Paths and cycles}

For paths $P_n$ ($n\ge 1$) and cycles $C_n$ ($n\ge 3$) the following non-trivial result holds:
$$\gamma_g(P_n) = \gamma_g(C_n) = \left\{
\begin{array}{ll}
\left\lceil\frac{n}{2}\right\rceil-1; & n \equiv 3\ (\bmod\ 4), \\
\\
\left\lceil\frac{n}{2}\right\rceil;& {\rm otherwise}.
\end{array}
\right.$$
The only published proof of  this theorem can be found in~\cite{kosmrlj-2017}.  The result implies that each of $P_n$ and $C_n$ is a $1/2$-graph if and only if $n\ (\bmod\ 4)\in \{0,1,2\}$.

\subsection*{Broken ladders}

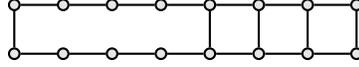
\begin{figure}[ht!]
	\begin{center}
		\begin{tikzpicture}[thick,scale=0.65]
		\tikzstyle{every node}=[circle, draw, fill=black!10,
		inner sep=0pt, minimum width=4pt]
		\begin{scope}
		\node (4) at (2,0) {};
		\node (5) at (2,1) {};
		\node (6) at (3,0) {};
		\node (7) at (3,1) {};
		\node (8) at (4,0) {};
		\node (9) at (4,1) {};
		\node (10) at (5,0) {};
		\node (11) at (5,1) {};
		\node (12) at (6,0) {};
		\node (13) at (6,1) {};
		\node (14) at (7,0) {};
		\node (15) at (7,1) {};
		\node (16) at (8,0) {};
		\node (17) at (8,1) {};
		\node (18) at (9,0) {};
		\node (19) at (9,1) {};
		\draw (4) -- (6) -- (8) -- (10) -- (12) -- (14) -- (16) -- (18) -- (19) -- (17) -- (15) -- (13) -- (11) -- (9) -- (7) -- (5) -- (4);
		\draw (17) -- (16);
		\draw (15) -- (14);
		\draw (13) -- (12);
		\end{scope}
		\end{tikzpicture}
		\caption{The broken ladder $BL_2$.}
		\label{fig:BL}
	\end{center}
\end{figure}

If $k\ge 1$, then the \emph{broken ladder} $BL_k$ is the graph obtained from the Cartesian product $P_4 \cp K_2$ by adding a path of length $4k+1$ between two adjacent vertices of degree $2$. See Fig.~\ref{fig:BL} for $BL_2$. Moreover, we set  $BL_0 = P_4 \cp K_2$.

\begin{proposition} {\rm \cite[Theorem 3.3]{kosmrlj-2014}}
	\label{prop:BL}
	If $k \geq 0$, then $\gg(BL_k) = 2 (k+2) = \frac{n(BL_k)}{2}$.
\end{proposition}


\subsection*{Hatted cycles}

\begin{figure}[ht!]
	\begin{center}
		\begin{tikzpicture}[thick,scale=0.65, rotate=90]
		
		\tikzstyle{every node}=[circle, draw, fill=black!10,
		inner sep=0pt, minimum width=4pt]
		
		\begin{scope}
		\node[label=below: $x'$] (0) at (0:1.5cm) {};
		\node[label=above: $x$] (1) at (0:2.5cm) {};
		\node[label=left: $y$] (2) at (40:2cm) {};
		\node (3) at (80:2cm) {};
		\node (4) at (120:2cm) {};
		\node (5) at (160:2cm) {};
		\node (6) at (200:2cm) {};
		\node (7) at (240:2cm) {};
		\node (8) at (280:2cm) {};
		\node[label=right: $y'$] (9) at (320:2cm) {};
		
		\draw (0) -- (2) -- (3) -- (4) -- (5) -- (6) -- (7) -- (8) -- (9) -- (0);
		\draw (9) -- (1) -- (2);
		\end{scope}
	
		\end{tikzpicture}
		\caption{The graph $\widehat{C}_9$.}
		\label{fig:cycle'}
	\end{center}
\end{figure}
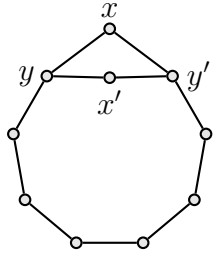

If $n\ge 4$, then the {\em hatted cycle} $\widehat{C}_n$ is obtained from the cycle $C_n$ by adding a new vertex and connecting it to two vertices at distance $2$ on the cycle; see Fig.~\ref{fig:cycle'} for $\widehat{C}_9$. (In~\cite{kosmrlj-2014} these graphs were denoted with $C_n'$.)

Using similar reasoning as in the proof of~\cite[Theorem 2.2]{kosmrlj-2014}, we now prove that hatted cycles of order $4k+2$ are $1/2$-graphs.

\begin{proposition}
	\label{prop:cycle-hat}
	If $k \geq 1$, then $\gg(\widehat{C}_{4k+1}) = 2k + 1 = \frac{n(\widehat{C}_{4k+1})}{2}$.
\end{proposition}

\proof
We use the notation from Fig.~\ref{fig:cycle'}. If Dominator starts the game on $d_1 = y$, then at most one of the vertices $x$ and $x'$ can be played in the remaining moves. Thus, after $d_1 = y$ is played, the game is the same as if it was played on $C_{4k+1}$. It follows that $\gg(\widehat{C}_{4k+1}) \leq \gg(C_{4k+1}) = 2k + 1$.

Before describing an optimal strategy of Staller, we define a {\em run} to be a maximal sequence of consecutive dominated vertices. Her strategy is to dominate only one new vertex in each of her moves. She achieves this by playing on the end of a run, or on $x$ or $x'$ if there is only one run and $y, y'$ are the end-vertices of this run. Let $m$ denote the number of moves in the D-game on $\widehat{C}_{4k+1}$.

If $m$ is even, then Staller dominates $\frac{m}{2}$ vertices and Dominator can dominate at most $4 + 3 (\frac{m}{2}-1)$ vertices. Together, both players dominate at most $2m + 1$ vertices, which must be at least $n(\widehat{C}_{4k+1}) = 4k+2$. Hence, $m \geq 2k + 1$.

If $m$ is odd, we use a similar reasoning as in the previous case. Together both players dominate at most $\frac{m-1}{2} + 4 + 3 (\frac{m+1}{2} - 1) = 2m + 2$ vertices. Hence, $m \geq 2k$. But as $m$ is odd, it cannot equal $2k$. Thus $m \geq 2k+1$ holds also in this case.

It follows that $\gg(\widehat{C}_{4k+1}) \geq 2k + 1$ and we conclude the equality.
\qed

The same reasoning proves that $\gg(\widehat{C}_{n}) = \left\lceil \frac{n(\widehat{C}_{n})}{2}\right\rceil -1$ if $n\ (\bmod\ 4)\in  \{0,2,3\}$. Therefore, every hatted cycle satisfies Conjecture~\ref{conj:Rall}, but only those which are obtained from $C_{4k+1}$ are $1/2$-graphs.

\section{Unicyclic traceable graphs}
\label{sec:unicyclic}

In this section we prove that Conjecture~\ref{conj:Rall} holds for all unicyclic traceable graphs.

Clearly, cycles are unicyclic traceable graphs.  On the other hand, if two nonadjacent vertices of the cycle of a unicyclic graph $G$ are of degree at least $3$, then $G$ is not traceable.  From this fact it is easy to deduce that if $G$ is traceable then $G$ is a cycle, or a graph obtained by attaching a path to a vertex of a cycle, or a graph obtained by attaching two disjoint paths to adjacent vertices of a cycle. The graphs from the latter two families will be called {\em tadpole graphs} and {\em two tailed tadpole graphs}, respective. Since we already know that Conjecture~\ref{conj:Rall} holds for cycles, to goal of this section is thus to prove that the conjecture also holds for tadpole graphs and for two tailed tadpole graphs.

From the preliminaries recall that the Union Lemma uses the weighting function $\w(P_n')$. To make computations simpler, we sometimes use  this function in the following equivalent form:
\[
\w(P'_{n})= \w(P''_{n})=
\begin{cases}
\frac{n}{2}; &  n \equiv 0\ (\bmod\ 4), \\
\frac{n}{2}+ \frac{1}{2}; &  n \equiv 1,2\ (\bmod\ 4),\\
\frac{n}{2}+\frac{1}{4}; &  n \equiv 3\ (\bmod\ 4).
\end{cases}
\]

\subsection*{Tadpole graphs}

If $m\ge 3$ and $n\ge 1$, then the \emph{$(m,n)$-tadpole graph} $T_{m,n}$ is the graph obtained from the disjoint union of a cycle $C_m$ and a path $P_n$ by joining a vertex of $C_m$ with an end-vertex of $P_n$. Clearly, $n(T_{m,n}) = n + m$.

\begin{theorem}
	\label{thm:tadpole}
	If $m \geq 3$ and $n \geq 1$, then $\gg(T_{m,n}) \leq \left\lceil\frac{m+n}{2}\right\rceil$.
\end{theorem}

\proof
Let $n = 4k + x + 1$ and $m = 4\ell + y + 3$, where $x, y \in \{0,1,2,3\}$. Let  $v$ be the vertex of $T_{m,n}$ of degree $3$. The first move $d_1 = v$ of Dominator implies that
$$\gg(T_{m,n}) \leq 1 + \gg'(P_{4k+x}' \cup P_{4\ell+y}'')\,.$$
Since $\gamma_g'(P_r') = \gamma_g'(P_r'')$ holds by Lemma~\ref{lem:P'-P''}, it suffices to consider the cases when $x\le y$, that is, ten such cases. Each of them can be handled using the Union Lemma.

Suppose first that $x = y = 0$. Then
$$\gg(T_{m,n}) \leq 1 + \gg'(P_{4k}' \cup P_{4\ell}'') \le 1 +\w(P_{4k}') + \w(P_{4\ell}'') = 1 + 2k + 2\ell\,,$$
where the second inequality follows by the Union Lemma. Since $n(T_{m,n}) = 4k + 4\ell + 4$ we get that $\gamma_g(T_{m,n}) \le \left\lceil\frac{m+n}{2}\right\rceil$.

Suppose next that $x = 1$ and $y = 3$. Then
\begin{align*}
\gg(T_{m,n}) & \leq 1 + \gg'(P_{4k+1}' \cup P_{4\ell+3}'') \le 1 +\w(P_{4k+1}') + \w(P_{4\ell+3}'') \\
& = 1 + (2k+1) + (2\ell+ 7/4) = 2k + 2\ell + 15/4\,.
\end{align*}
Since $n(T_{m,n}) = 4k + 4\ell + 8$ we get that $\gamma_g(T_{m,n}) \le \left\lceil\frac{m+n}{2}\right\rceil$.

\begin{table}[ht!]
	\begin{center}
		\begin{tabular}{c|c|c|c|c}
			$x$ & $y$ & $1 + \lceil w(P'_{4k + x}) + w(P''_{4\ell + y}) \rceil$ & $n(T_{m,n})$ & $\lceil \frac{n(T_{m,n})}{2} \rceil$ \\ \hline
			0 & 0 & $\lceil 2 k+2 \ell\rceil +1 = 2k + 2\ell + 1$ & $4k + 4\ell + 4$ & $2k + 2\ell + 2$ \\
			0 & 1 & $\lceil 2 k+2 \ell\rceil +2 = 2k + 2\ell + 2$ & $4k + 4\ell + 5$ & $2k + 2\ell + 3$ \\
			0 & 2 & $\left\lceil 2 k+2 \ell+\frac{3}{2}\right\rceil +1 = 2k + 2\ell + 3$ & $4k + 4\ell + 6$ & $2k + 2\ell + 3$ \\
			0 & 3 & $\left\lceil 2 k+2 \ell+\frac{7}{4}\right\rceil +1 = 2k + 2\ell + 3$ & $4k + 4\ell + 7$ & $2k + 2\ell + 4$ \\
			1 & 0 & $\lceil 2 k+2 \ell \rceil +2 = 2k + 2\ell + 2$ & $4k + 4\ell + 5$ & $2k + 2\ell + 3$ \\
			1 & 1 & $\lceil 2 k+2 \ell \rceil +3 = 2k + 2\ell + 3$ & $4k + 4\ell + 6$ & $2k + 2\ell + 3$ \\
			1 & 2 & $\left\lceil 2 k+2 \ell+\frac{5}{2}\right\rceil +1 = 2k + 2\ell + 4$ & $4k + 4\ell + 7$ & $2k + 2\ell + 4$ \\
			1 & 3 & $\left\lceil 2 k+2 \ell+\frac{11}{4}\right\rceil +1 = 2k + 2\ell + 4$ & $4k + 4\ell + 8$ & $2k + 2\ell + 4$ \\
			2 & 0 & $\left\lceil 2 k+2 \ell+\frac{3}{2}\right\rceil +1 = 2k + 2\ell + 3$ & $4k + 4\ell + 6$ & $2k + 2\ell + 3$ \\
			2 & 1 & $\left\lceil 2 k+2 \ell+\frac{5}{2}\right\rceil +1 = 2k + 2\ell + 4$ & $4k + 4\ell + 7$ & $2k + 2\ell + 4$ \\
			2 & 2 & $\lceil 2 k+2 \ell\rceil +4 = 2k + 2\ell + 4$ & $4k + 4\ell + 8$ & $2k + 2\ell + 4$ \\
			2 & 3 & $\left\lceil 2 k+2 \ell+\frac{13}{4}\right\rceil +1 = 2k + 2\ell + 5$ & $4k + 4\ell + 9$ & $2k + 2\ell + 5$ \\
			3 & 0 & $\left\lceil 2 k+2 \ell+\frac{7}{4}\right\rceil +1 = 2k + 2\ell + 3$ & $4k + 4\ell + 7$ & $2k + 2\ell + 4$ \\
			3 & 1 & $\left\lceil 2 k+2 \ell+\frac{11}{4}\right\rceil +1 = 2k + 2\ell + 4$ & $4k + 4\ell + 8$ & $2k + 2\ell + 4$ \\
			3 & 2 & $\left\lceil 2 k+2 \ell+\frac{13}{4}\right\rceil +1 = 2k + 2\ell + 5$ & $4k + 4\ell + 9$ & $2k + 2\ell + 5$ \\
			3 & 3 & $\left\lceil 2 k+2 \ell+\frac{7}{2}\right\rceil +1 = 2k + 2\ell + 5$ & $4k + 4\ell + 10$ & $2k + 2\ell + 5$ \\
		\end{tabular}
		\caption{The calculations for all different cases.}
		\label{table:tadpole}
	\end{center}
\end{table}

The remaining cases to be considered can be treated along the same lines as the above two cases. In Table~\ref{table:tadpole}  the summary of calculations for all the cases is presented.
\qed

\subsection*{Two tailed tadpole graphs}

If $m\ge 3$ and $n, k\ge 1$, then the notation $T_{m,n,k}$ means a graph obtained from the disjoint union of a cycle $C_m$ and paths $P_n$ and $P_k$ by joining adjacent vertices of $C_m$ with end-vertices of $P_n$ and $P_k$ with two independent edges. Clearly, $n(T_{m,n,k}) = n + m + k$.

\begin{theorem}
	\label{thm:2tailed-tadpole}
	If $m \geq 3$ and $n,k \geq 1$, then $\gg(T_{m,n,k}) \leq \left\lceil\frac{n(T_{m,n,k})}{2}\right\rceil$.
\end{theorem}

\proof
Let $G = T_{m,n,k}$ and let the vertices of $G$ be denoted by $v_1, \dots, v_{n+m+k}$ such that $v_1\dots v_{n+m+k}$ is the Hamiltonian path and $v_{n+1}v_{n+m}$ is the extra edge. During the game, an \emph{antirun} is a component of the subgraph induced by the undominated vertices.

Suppose that $d_1 = v_{n+1}$. After this move, we have three antiruns: $X=v_1\dots v_{n-1}$,\, $Y= v_{n+3}\dots v_{n+m-1}$,\, and $Z=v_{n+m+1}\dots v_{n+m+k}$. Together with the neighboring dominated vertices, we may consider the antiruns $X$, $Y$, $Z$ as a $P_{n-1}'$, a $P_{m-3}''$, and a $P_{k}'$, respectively. By the Continuation Principle, we can consider antiruns $Y$ and $Z$ as one, and get $\gg'(G|N[d_1]) \leq \gg'(P'_{n-1} \cup P'_{m+k-2})$. Denote $n = 4n'+x+1$, $k=4k'+z$, $m = 4m' + y + 2$, where $x,y,z \in \{0,1,2,3\}$ and $m',n',k'$ are integers. Hence the Union Lemma yields
 $$\gg(G) \leq 1 + \gg'(P'_{n-1} \cup P'_{m+k-2}) \leq 1 + \left\lceil \w(P'_{4n'+x}) + \w(P'_{4(m'+k')+y+z}) \right\rceil \,.$$

Repeating similar calculations as in the proof of Theorem~\ref{thm:tadpole} for all $64$ different values of $(x,y,z)$, we see that for most cases $\gg(G) \leq \left \lceil \frac{n(G)}{2} \right \rceil$ holds. The exceptional cases, after transforming from $(x,y,z)$ to $(n,m,k)\ (\bmod\ 4)$, are gathered in Table~\ref{table:exceptions}.

\begin{table}[ht!]
	\begin{center}
		\begin{tabular}{r|cccccccccccccccc}
			$n\ (\bmod\ 4)$ &2&2&2&2&3&3&3&3&3&3&3&3&0&0&0&0 \\ \hline
			$m\ (\bmod\ 4)$ &2&3&0&1&2&2&3&3&0&0&1&1&2&3&0&1 \\ \hline
			$k\ (\bmod\ 4)$ &2&1&0&3&1&3&0&2&1&3&0&2&2&1&0&3 \\
		\end{tabular}
		\caption{Exceptional cases.}
		\label{table:exceptions}
	\end{center}
\end{table}

Most of those cases can be omitted by symmetry (assuming $d_1 = v_{n+m}$ and then repeating the calculations using the Union Lemma). The only problematic cases left are:

\begin{itemize}
	\itemsep0em
	\item $m \equiv 2$, $n \equiv 2$, and $k \equiv 2\ (\bmod\ 4)$,
	\item $m \equiv 1$, $n \equiv 2$, and $k \equiv 3\ (\bmod\ 4)$,
	\item $m \equiv 2$, $n \equiv 3$, and $k \equiv 3\ (\bmod\ 4)$,
	\item $m \equiv 0$, $n \equiv 3$, and $k \equiv 3\ (\bmod\ 4)$,
	\item $m \equiv 1$, $n \equiv 3$, and $k \equiv 0\ (\bmod\ 4)$,
	\item $m \equiv 0$, $n \equiv 0$, and $k \equiv 0\ (\bmod\ 4)$,
\end{itemize}

First consider the case $m \equiv 2$, $n \equiv 2$, and $k \equiv 2\ (\bmod\ 4)$. Suppose Dominator starts the game on $d_1 = v_n$. His second move $d_2$ is such that after $s_1$ and $d_2$, all three vertices in $\{ v_{n+m}, v_{n+m+1}, v_{n+m+2} \}$ are dominated. By using the Continuation Principle, the orders of the antiruns are $n-2 \equiv 0$, $m-2 \equiv 0$, and $k-2 \equiv 0\ (\bmod\ 4)$. The Union Lemma then yields
$$\gg(G) \leq 3 + \left \lceil \frac{n-2}{2} + \frac{m-2}{2} + \frac{k-2}{2} \right \rceil = 3 + \left \lceil \frac{n+m+k}{2} -3 \right \rceil = \frac{n(G)}{2}\,.$$

For all the remaining cases, we assume again that $d_1=v_{n+1}$ is the first move of Dominator. His second move $d_2$ will be specified depending on the move $s_1$ of Staller.

At the beginning, there is only one antirun. After the move $d_1=v_{n+1}$ that dominates four vertices, we have the following three antiruns: $X=v_1\dots v_{n-1}$,\, $Y= v_{n+3}\dots v_{n+m-1}$,\, and $Z=v_{n+m+1}\dots v_{n+m+k}$. Together with the neighboring dominated vertices, we may consider the antiruns $X$, $Y$, $Z$ as a $P_{n-1}'$, a $P_{m-3}''$, and a $P_{k}'$, respectively. Note that $X$ or $Y$ might be a path of order $0$ (i.e., an empty graph), but we always assume that $n,m,k$ are positive integers. After the move $d_1$, at any point in the game, each antirun will be a path $P_j'$ or $P_j''$ for an appropriate $j$.

In all the remaining cases $n(G)$ is even. Suppose that Staller's move $s_1$  increases the number of antiruns. Then, we have four antiruns $X_1, \dots , X_4$ and, as such a move $s_1$ increases the number of dominated vertices by $3$, we have $|X_1|+|X_2|+|X_3|+|X_4|= n+m+k-7$. We consider three cases.  If there is an antirun containing at least three vertices, Dominator can choose $d_2$ such that the number of antiruns is not increased and the move dominates three new vertices. Applying the Union Lemma for this graph with antiruns $X_1^*, \dots , X_4^*$, we obtain the following estimation on the total number of moves $t$:
\begin{align*}
t  & \le 3+ \gamma_g'\left(\bigcup_{i = 1}^4 X_i^* \right) \le 3+  \left\lceil  \sum_{i = 1}^4 \w(X_i) \right\rceil \le 3+  \left\lceil  \sum_{i = 1}^4 \left(\frac{|X_i^*|}{2}+\frac{1}{2}\right) \right\rceil\\
& = 3+ \left\lceil \frac{n+m+k-10}{2}+4\cdot \frac{1}{2} \right\rceil
=3+ \left\lceil \frac{n(G)}{2}-3 \right\rceil = \frac{n(G)}{2}\,.
\end{align*}
Note that the last step uses the fact that $n(G)$ is even. In the second case, after Staller's move $s_1$, there is no antirun of order at least three but there is an antirun with $|X_i|=2$. Then, Dominator may play a vertex $d_2$ that dominates  the entire $X_i$ and, for the remaining antiruns $X_1^*$, $X_2^*$, and $X_3^*$, we have
\begin{align*}
t  & \le 3+ \gamma_g'\left(\bigcup_{i = 1}^3 X_i^* \right)  \le 3 +  \left\lceil  \sum_{i = 1}^3 \left(\frac{|X_i^*|}{2}+\frac{1}{2}\right) \right\rceil \\
& = 3+ \left\lceil \frac{n+m+k-9}{2}+3\cdot \frac{1}{2} \right\rceil = \frac{n(G)}{2}\,.
\end{align*}
In the third case, every antirun consists of one vertex after Staller's move $s_1$ and, therefore, $n(G)=11$. It can be checked by hand (we also checked by computer) that the game can be finished in $6$ moves, thus $\gg(G) \le \lceil n(G)/2 \rceil$ holds.
\medskip

From now on, we may assume that Staller's move $s_1$ does not increase the number of antiruns. Let $X' \subseteq X$,\, $Y' \subseteq Y$, and $Z' \subseteq Z$ be the antiruns after the move $s_1$. Remark that $n(G)=n+m+k$ is even for each of the following cases.
\begin{itemize}
	\item \textbf{Case 1.} $m \equiv 1$, $n \equiv 2$, and $k \equiv 3\ (\bmod\ 4)$ \\
	If Staller dominates at least one vertex from $X$, we may assume by the Continuation Principle that $|X'|=|X|-1= n-2$ and then $|X'|\equiv 0 (\bmod\ 4)$. Then, Dominator responds by playing $d_2=v_{n+m+2}$ which dominates three vertices from $Z'=Z$. The order of the antiruns are
	$n-2 \equiv 0\,$,
	$m-3 \equiv 2\, $,
	$k-3 \equiv 0 (\bmod\ 4)$ \,
	and, by the Union Lemma, we get
	\begin{align*}
	t & \le  3+  \left\lceil  \frac{n-2}{2} + \frac{m-3}{2}+ \frac{1}{2}+\frac{k-3}{2} \right\rceil= 3+ \left\lceil \frac{n+m+k}{2} - \frac{7}{2} \right\rceil = \frac{n(G)}{2}.
	\end{align*}
	If Staller dominates at least one vertex from $Y$, we may apply the Continuation Principle again. Then, Dominator  plays $d_2=v_{n+m+2}$ which dominates three vertices. The antiruns are of the following orders:
	$n-1 \equiv 1\,$,
	$m-4 \equiv 1\, $,
	$k-3 \equiv 0\ (\bmod\ 4)$.
	By the Union Lemma,
	$$
	t  \le  3+  \left\lceil  \frac{n-1}{2} +\frac{1}{2} + \frac{m-4}{2}+ \frac{1}{2}+\frac{k-3}{2} \right\rceil= 3+ \left\lceil \frac{n+m+k}{2} - 3 \right\rceil = \frac{n(G)}{2}.	
	$$
	In the third case, Staller does not dominate any vertices from $X \cup Y$ and consequently, she  has to dominate at least two vertices from $Z$ by playing either $v_{m+n+1}$, $v_{m+n+2}$, $v_{m+n+k-1}$, or $v_{m+n+k}$. By the Continuation Principle, we may assume that $|Z'|=|Z|-2$. In the next move, Dominator plays $d_2=v_{n+2}$ and dominates two vertices from $Y'=Y$.  This creates antiruns with the following orders:
	$n-1 \equiv 1\,$,
	$m-5 \equiv 0\, $,
	$k-2 \equiv 1\ (\bmod\ 4)$. By the Union Lemma,
	$$
	t  \le  3+  \left\lceil  \frac{n-1}{2} +\frac{1}{2} + \frac{m-5}{2}+ \frac{k-2}{2} +\frac{1}{2} \right\rceil= 3+ \left\lceil \frac{n+m+k}{2} - 3 \right\rceil = \frac{n(G)}{2}.	
	$$
	This finishes the proof for Case 1, since it follows that $\gg(G) \leq t \leq \frac{n(G)}{2}$.
	
	\item \textbf{Case 2.} $m \equiv 2$, $n \equiv 3$, and $k \equiv 3\ (\bmod\ 4)$ \\
	Similarly to the previous case, we consider three subcases according to Staller's move $s_1$. If Staller dominates at least one vertex from $X$, Dominator replies with  $d_2=v_{n+m+2}$. After this move, the orders of the antiruns are
	$n-2 \equiv 1\,$,
	$m-3 \equiv 3\, $, and
	$k-3 \equiv 0\ (\bmod\ 4)$.
	The counting gives
	$$
	t  \le  3+  \left\lceil  \frac{n-2}{2} +\frac{1}{2} + \frac{m-3}{2}+ \frac{1}{4} + \frac{k-3}{2}  \right\rceil= 3+ \left\lceil \frac{n+m+k}{2} - \frac{13}{4} \right\rceil = \frac{n(G)}{2}	
	$$
	by the Union Lemma.
	If Staller dominates at least one vertex from $Y$, Dominator plays $d_2=v_{n+m+2}$ again and we have antiruns satisfying
	$n-1 \equiv 2\,$,
	$m-4 \equiv 2\, $, and
	$k-3 \equiv 0\ (\bmod\ 4)$.
	The inequality gives
	$$
	t  \le  3+  \left\lceil  \frac{n-1}{2} +\frac{1}{2} + \frac{m-4}{2}+ \frac{1}{2} + \frac{k-3}{2}  \right\rceil= 3+ \left\lceil \frac{n+m+k}{2} - 3 \right\rceil = \frac{n(G)}{2}.	
	$$
	If Staller does not dominate any vertices from $X \cup Y$, then she  dominates at least two vertices from $Z$. We may assume by the Continuation Principle that $|Z'|=|Z|-2$. Then, Dominator plays $d_2=v_{n+3}$ and dominates three vertices from $Y'=Y$. After this move we have the following antiruns:
	$n-1 \equiv 2\,$,
	$m-6 \equiv 0\, $,
	$k-2 \equiv 1\ (\bmod\ 4)$. The counting  gives
	$$
	t \le  3+  \left\lceil  \frac{n-1}{2} +\frac{1}{2} + \frac{m-6}{2}+  \frac{k-2}{2} +\frac{1}{2}  \right\rceil= 3+ \left\lceil \frac{n+m+k}{2} - \frac{7}{2} \right\rceil = \frac{n(G)}{2}.	
	$$
	
	\item \textbf{Case 3.} $m \equiv 0$, $n \equiv 3$, and $k \equiv 3\ (\bmod\ 4)$ \\
	As the proof for the remaining cases are very similar to the previous ones, we describe only the main points from the argument.
	If Staller dominates at least one vertex from $X$, Dominator plays $d_2=v_{n+m+2}$ and then, the antiruns are
	$n-2 \equiv 1\,$,
	$m-3 \equiv 1\, $, and
	$k-3 \equiv 0\ (\bmod\ 4)$. If Staller dominates at least one vertex from $Y$, Dominator plays $d_2=v_{n+m+2}$ again. Then, we have antiruns of order
	$n-1 \equiv 2\,$,
	$m-4 \equiv 0\, $, and
	$k-3 \equiv 0\ (\bmod\ 4)$. If Staller plays a vertex from $Z$, Dominator replies by dominating at least two vertices from $X$. The antiruns are:
	$n-3 \equiv 0\,$,
	$m-3 \equiv 1\, $, and
	$k-2 \equiv 1\ (\bmod\ 4)$. By applying the Union Lemma in the continuation, we obtain $t \le n(G)/2$
	for all cases.
	
	\item \textbf{Case 4.} $m \equiv 1$, $n \equiv 0$, and $k \equiv 3\ (\bmod\ 4)$ \\
	If Staller plays $s_1$ from $X$, Dominator selects $d_2=v_{n+m+2}$ and then, we have the antiruns
	$n-2 \equiv 2$,
	$m-3 \equiv 2$, and
	$k-3 \equiv 0\ (\bmod\ 4)$.
	If $s_1$ dominates at least one vertex from $Y$, Dominator plays	$d_2=v_{n+m+2}$ and then, the antiruns are
	$n-1 \equiv 3$,
	$m-4 \equiv 1$, and
	$k-3 \equiv 0\ (\bmod\ 4)$.
	Finally, if Staller selects a vertex $s_1$ and dominates at least two vertices from $Z$, Dominator can play $d_2=v_2$. Since $n >0$ and $n \equiv 0\ (\bmod\ 4)$, he dominates at least three vertices with this move. The antiruns are of size $n-4 \equiv 0$,
	$m-3 \equiv 2$, and
	$k-2 \equiv 1\ (\bmod\ 4)$. The desired inequality $t \le n(G)/2$ follows.
	
	\item \textbf{Case 5.} $m \equiv 0$, $n \equiv 0$, and $k \equiv 0\ (\bmod\ 4)$ \\
	We may assume without loss of generality that $n \ge k$. (Otherwise, Dominator would play $v_{m+n}$ instead of $v_{n+1}$.) If Staller's move dominates a vertex from $X$ and $X'$ still contains at least three undominated vertices, then Dominator plays in $X'$ and leaves $|X|-5$ undominated vertices there. It follows that the antiruns are of order
	$n-5 \equiv 3$,
	$m-3 \equiv 1$, and
	$k \equiv 0\ (\bmod\ 4)$. The usual counting then proves that the number of moves is at most $n(G)/2$. If Dominator cannot do this, then $n-2=2$ and he tries to dominate three vertices from $Y'$ if possible. In this case, the orders of the antiruns are  $n-2 = 2$,
	$m-6 \equiv 2$,
	$k \equiv 0\ (\bmod\ 4)$,
	and the result follows as earlier. If neither  $X'$ nor $Y'$ contains three undominated vertices, then $n=4$, $m=4$ and, by the assumption $n \ge k >0$ and $k \equiv 0\ (\bmod\ 4)$, $k=4$ follows. It can be checked by hand (or computer), that the game domination number of this graph is $6=n(G)/2$. For the remaining cases Dominator's startegy is the usual. If Staller dominates a vertex from $Y$, then Dominator dominates three vertices from $X'$. Note that this is possible as the antirun contains $n-1 \equiv 3\ (\bmod\ 4)$ vertices after Staller's move. If Staller dominates at least two vertices from $Z$, Dominator can reply by dominating at least three vertices from $X'$.
\end{itemize}
Since for each possible case, we have shown a strategy of Dominator which ensures that the game is finished in at most $\lceil n(G)/2 \rceil$ moves, the desired upper bound $\gg(G) \le \lceil n(G)/2 \rceil$ follows.  \qed

\section{Additional classes supporting Conjecture~\ref{conj:Rall}}
\label{sec:additional}

Applying a relationship between the domination game and minimal edge cuts, non-trivial families of $1/2$-graphs were constructed in~\cite[Section~5]{klavzar-2019} that support Conjecture~\ref{conj:Rall}. In this section we find additional non-trivial families that support the conjecture.

\subsection*{Cycles with a chord}

A graph obtained from a cycle $C_n$ by adding an edge between two nonadjacent vertices is clearly traceable. These graphs form our next family for which Conjecture~\ref{conj:Rall} holds.

\begin{proposition}
	\label{prop:cycles-chord}
	If $G$ is a graph obtained from a cycle $C_n$ by connecting two nonadjacent vertices of the cycle, then $\gg(G) \leq \left\lceil \frac{n(G)}{2}\right\rceil$.
\end{proposition}

\proof
Let $v_1, \ldots, v_n$ be the vertices of the cycle and let $v_1 v_i$, $3 \leq i \leq n-1$, be the additional edge. Suppose Dominator starts the game by playing $d_1 = v_1$. Then it follows by the Continuation Principle that
$$\gg(G) \leq 1 + \gg'(P''_{n-3}|v_i) \leq 1 + \gg(P''_{n-3})\,.$$
Applying Lemma~\ref{lem:P'-P''} in each of the following cases, we get:
\begin{itemize}
	\item If $n \equiv 0 \text{ or } 2\ (\bmod\ 4)$, then $$\gg(G) \leq 1 + \left \lceil \frac{n-3}{2} \right \rceil = 1 + \frac{n-2}{2} = \frac{n}{2} = \left \lceil \frac{n}{2} \right \rceil\,.$$
	
	\item If $n \equiv 3\ (\bmod\ 4)$, then $$\gg(G) \leq 1 + \left \lceil \frac{n-3}{2} \right \rceil = 1 + \frac{n-3}{2} = \frac{n-1}{2} \leq \left \lceil \frac{n}{2} \right \rceil\,.$$
	
	\item If $n \equiv 1\ (\bmod\ 4)$, then $$\gg(G) \leq 1 + \left \lceil \frac{n-3}{2} \right \rceil + 1 = 2 + \frac{n-3}{2} = \frac{n+1}{2} = \left \lceil \frac{n}{2} \right \rceil\,.$$
\end{itemize}
\qed

\subsection*{Graphs from ${\cal F}(X)$, where $X$ is traceable}

\begin{figure}[ht!]
	\begin{center}
		\begin{tikzpicture}[thick,scale=0.65, rotate=90]
		
		\tikzstyle{every node}=[circle, draw, fill=black!10,
		inner sep=0pt, minimum width=4pt]
		
			\begin{scope}
			\node (0) at (0:1cm) {};
			\node (1) at (0:2cm) {};
			\node (00) at (0:3cm) {};
			\node (11) at (0:4cm) {};
			\node[label=left: $y$] (2) at (40:2cm) {};
			\node (3) at (80:2cm) {};
			\node (4) at (120:2cm) {};
			\node (5) at (160:2cm) {};
			\node (6) at (200:2cm) {};
			\node (7) at (240:2cm) {};
			\node (8) at (280:2cm) {};
			\node[label=right: $y'$] (9) at (320:2cm) {};
			
			\draw (0) -- (2) -- (3) -- (4) -- (5) -- (6) -- (7) -- (8) -- (9) -- (0);
			\draw (9) -- (1) -- (2);
			\draw (9) -- (00) -- (2);
			\draw (9) -- (11) -- (2);
			\draw (0) -- (1) -- (00) -- (11);
			\end{scope}
		
		\end{tikzpicture}
		\caption{A graph from ${\cal F}(P_4)$.}
		\label{fig:cycle'2}
	\end{center}
\end{figure}
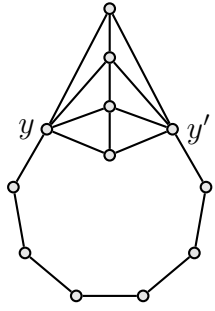

Let $X$ be a traceable graph. Then the family ${\cal F}(X)$ consists of all graphs $G$ that can be constructed in the following way. $G$ is obtained from the disjoint union of $X$ and a path $P_{n}$, $n\ge 3$, with end-vertices $y$ and $y'$, by connecting $y$ to all vertices of $X$, and by connecting $y'$ to the vertices from $W\subseteq V(X)$, where $W$ contains at least one end-vertex of some Hamiltonian path in $X$. In the example in Fig.~\ref{fig:cycle'2} we have $X = P_4$, $n = 8$,  and $y'$ is adjacent to all vertices of $P_4$, that is, $W = V(P_4)$.

If $X$ is a traceable graph and $G\in {\cal F}(X)$, then it is easy to observe that $G$ is traceable. Hence the  following result is of interest to us.

\begin{proposition}
	\label{prop:generalized-hat}
If $X$ is a traceable graph and $G\in {\cal F}(X)$, then $\gg(G) \leq \left\lceil \frac{n(G)}{2}\right\rceil$.
\end{proposition}

\proof
Let $X$ be a traceable graph and $G\in {\cal F}(X)$, where the end-vertex $y$ of the building graph $P_{n}$ of $G$ is adjacent to all the vertices of the building graph $X$ of $G$. Consider the D-game and set $d_1 = y$. This move dominates all the vertices of $X$ in $G$, and at most one of the vertices from $V(X)$ can be played in the rest of the game. If follows that after the first move $d_1 = y$ is played, the game is the same as if it would be played on $C_{n+1}$.  Since cycles fulfill Conjecture~\ref{conj:Rall}, the same holds for $G$.
\qed

\subsection*{Particular Halin graphs}

A {\em Halin graph} is a graph obtained from a plane embedding of a tree $T$ on at least four vertices and with no vertex of degree $2$, by connecting the leaves of $T$ into a cycle in the clock-wise order with respect to the embedding. These graphs were introduced in~\cite{halin-1971} and are of continuing interest, cf.~\cite{chen-2017}. For us the most important property of these graphs is that they are Hamiltonian~\cite{bondy-1985}.

\begin{figure}[ht!]
	\begin{center}
		\begin{tikzpicture}[thick,scale=0.27]
		\tikzstyle{every node}=[circle, draw, fill=black!10,
		inner sep=0pt, minimum width=3pt]
		
		\pgfmathtruncatemacro{\a}{24}
		\pgfmathtruncatemacro{\b}{8}
		\pgfmathtruncatemacro{\c}{4}
		\pgfmathtruncatemacro{\d}{3}
		
		\begin{scope}
		\node[label=above:$r$] (r) at (12.5, 3*\d) {};
		
		\foreach \x in {1,...,\a}
		\node (\x) at (\x, 0) {};
		
		\foreach \x in {1,...,\b}
		\node (\x11) at (3*\x - 1, \d) {};
		
		\foreach \x in {1,...,\c}
		\node (\x22) at (6*\x - 2.5, 2*\d) {};
		
		\foreach \x in {1,...,\c}
		\draw (\x22) -- (r);
		
		\draw (1) -- (111) -- (122) -- (211) -- (6);
		\draw (7) -- (311) -- (222) -- (411) -- (12);
		\draw (13) -- (511) -- (322) -- (611) -- (18);
		\draw (19) -- (711) -- (422) -- (811) -- (24);
		
		\draw (2) -- (111) -- (3);
		\draw (4) -- (211) -- (5);
		\draw (8) -- (311) -- (9);
		\draw (10) -- (411) -- (11);
		\draw (14) -- (511) -- (15);
		\draw (16) -- (611) -- (17);
		\draw (20) -- (711) -- (21);
		\draw (22) -- (811) -- (23);
		\end{scope}
		
		\begin{scope}[xshift=28cm]
		\node[label=above:$r$] (r) at (12.5, 3*\d) {};
		
		\foreach \x in {1,...,\a}
		\node (\x) at (\x, 0) {};

	    \foreach \x in {1,...,\b}
		\node (\x11) at (3*\x - 1, \d) {};
		
		\foreach \x in {1,...,\c}
		\node (\x22) at (6*\x - 2.5, 2*\d) {};
		
		\foreach \x in {1,...,\c}
		\draw (\x22) -- (r);
		
		\draw (1) -- (111) -- (122) -- (211) -- (6);
		\draw (7) -- (311) -- (222) -- (411) -- (12);
		\draw (13) -- (511) -- (322) -- (611) -- (18);
		\draw (19) -- (711) -- (422) -- (811) -- (24);
		
		\draw (2) -- (111) -- (3);
		\draw (4) -- (211) -- (5);
		\draw (8) -- (311) -- (9);
		\draw (10) -- (411) -- (11);
		\draw (14) -- (511) -- (15);
		\draw (16) -- (611) -- (17);
		\draw (20) -- (711) -- (21);
		\draw (22) -- (811) -- (23);

         \draw (1) -- (24);
         		\foreach \x in {1,...,\a}
		\node (\x) at (\x, 0) {};
		
		\draw (1) to[out=-35,in=-145] (\a);
		
		\end{scope}
		
		\end{tikzpicture}
		\caption{The graphs $T(3; 4,2,3)$ (left) and $H(3; 4,2,3)$ (right).}
		\label{fig:tree-and-Halin}
	\end{center}
\end{figure}
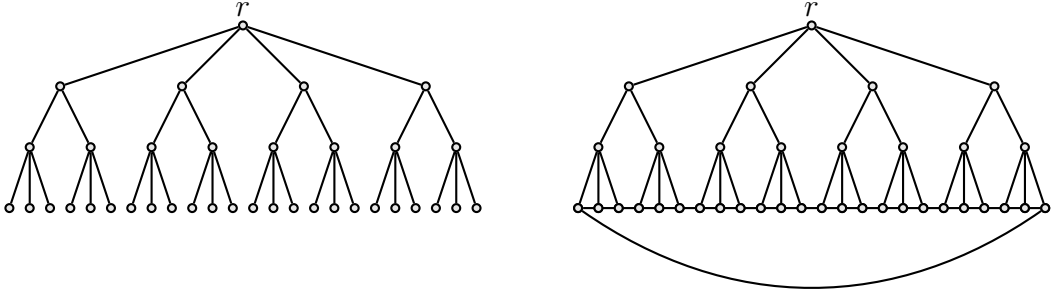

Let $k\ge 1$, $d_0\ge 3$, and $d_i\ge 2$ for $i\in [k-1]$. Then let $H(k; d_0, \ldots, d_{k-1})$ be the Halin graph obtained from the tree $T = T(k; d_0, \ldots, d_{k-1})$ defined as follows. Let $r$ be the root of $T$ of degree $d_0$. For $i\in [k-1]$, each vertex at distance $i$ from $r$ is of degree $d_i+1$. The vertices at distance $k$ from $r$ are the leaves of $T$. See Fig.~\ref{fig:tree-and-Halin} for $T(3; 4,2,3)$ and $H(3; 4,2,3)$.

\begin{proposition}
\label{prop:gamma-of-halin}
If $k\ge 1$, $d_i\ge 3$ for $i\in [k]_0$, and $H = H(k; d_0, \ldots, d_{k-1})$, then $\gamma(H) < \frac{n(H)}{4}$. Consequently, $\gamma_g(H) < \frac{n(H)}{2} - 1$.
\end{proposition}

\proof
Let $V_i$, $i\in [k+1]_0$, be the set of vertices of $H$ at distance $i$ from the root $r$. In particular, $V_0 = \{r\}$. Let
$$ D =
\begin{cases}
\vspace*{2mm}
V_0 \cup \bigcup_{i=1}^{k/3} V_{3i-1}\,; &  k \equiv 0\ (\bmod\ 3),\\
\vspace*{2mm}
V_0 \cup \bigcup_{i=1}^{(k-1)/3} V_{3i}\,; & k \equiv 1\ (\bmod\ 3),\\
\bigcup_{i=1}^{(k+1)/3} V_{3i-2}\,; & k \equiv 2\ (\bmod\ 3).
\end{cases}
$$
Since for every $i\in [k-1]$, the set $V_i$ dominates $V_{i-1}\cup V_i \cup V_{i+1}$, it readily follows that $D$ is a dominating set of $H$. Moreover, as $d_i\ge 3$, we also have that $|V_i| < |V_{i-1}\cup V_i \cup V_{i+1}|/4$, from which the first assertion of the proposition follows.

Since for every graph $G$ we have $\gamma_g(G) \le 2\gamma(G) - 1$ (see~\cite[Theorem 1]{bresar-2010}), the game domination number of $H$ can be bounded   as follows: $\gamma_g(H) \le 2\gamma(H) - 1 < 2\cdot \frac{n(H)}{4} - 1 = \frac{n(H)}{2} - 1$.
\qed

Note that the only requirement that the proof of Proposition~\ref{prop:gamma-of-halin} works is that  $|V_i| < |V_{i-1}\cup V_i \cup V_{i+1}|/4$ holds for those indices $i$ for which $V_i\subseteq D$. This is achieved by only requiring that $d_i\ge 3$ for the corresponding $i$.

\section{Computer support and a phenomenon}
\label{sec:computer}

A possible approach to prove (or disprove) Rall's conjecture is the following. Let $G$ be an arbitrary traceable graph and let $P$ be a Hamiltonian path in $G$. Then we know that the conjecture holds when the game is played on $P$. Adding edges to $P$ one by one, until $G$ is reached, while keeping the game domination number below $\lceil n(G)/2\rceil$, would yield the conjecture. An obstruction with this approach is the fact~\cite[Proposition 2.4]{bresar-2014} which asserts that for any $\ell\ge 5$ there exists a graph $G$ with an edge $e$ such that $\gamma_g(G) = \ell$ and $\gamma_g(G-e) = \ell - 2$. Thus, going from $G-e$ to $G$, the result implies that adding an edge to a graph, the game domination number can increase by up to $2$. But is the situation different for traceable graphs?  Moreover, we still have all possible Hamiltonian paths to start with, as well as all possible orders of the edges not on $P$ to be put back into $G$. In this section we report our computer experiments on this approach.

Our first partial support for Rall's conjecture was obtained by computer.

\begin{proposition}
	\label{prop:computer_paths}
	If $4 \leq n \leq 21$ and $G$ is a path $P_n$ with two additional edges, then $\gamma_g(G) \leq \lceil \frac{n}{2} \rceil$. If $4 \leq n \leq 15$, the same holds for a path $P_n$ with three additional edges.
\end{proposition}

As mentioned above, proving that in the case of traceable graphs adding edges to the graph does not increase the game domination number would suffice to prove Rall's conjecture. However, the following example shows that this is not true.

\begin{figure}[ht!]
	\begin{center}
		\begin{tikzpicture}[thick,scale=0.6]
		
		\tikzstyle{every node}=[circle, draw, fill=black!10,
		inner sep=0pt, minimum width=4pt]
		
		\begin{scope}
		\node (0) at (0,0) {};
		\node (1) at (1,0) {};
		\node (2) at (2,0) {};
		\node (3) at (3,0) {};
		\node (4) at (4,0) {};
		\node (5) at (5,0) {};
		\node (6) at (6,0) {};
		\node (7) at (7,0) {};
		\node (8) at (8,0) {};
		\node (9) at (9,0) {};
		\node (10) at (10,0) {};
		
		\draw (0) -- (1) -- (2) -- (3) -- (4) -- (5) -- (6) -- (7) -- (8) -- (9) -- (10);
		\draw (0) to[out=35,in=145] (4);
		\draw (5) to[out=35,in=145] (8);
		\draw (1) to[out=-35,in=-145] (7);
		\end{scope}
		
		\begin{scope}[xshift=12cm]
		\node (0) at (0,0) {};
		\node (1) at (1,0) {};
		\node (2) at (2,0) {};
		\node (3) at (3,0) {};
		\node (4) at (4,0) {};
		\node (5) at (5,0) {};
		\node (6) at (6,0) {};
		\node (7) at (7,0) {};
		\node (8) at (8,0) {};
		\node (9) at (9,0) {};
		\node (10) at (10,0) {};
		
		\draw (0) -- (1) -- (2) -- (3) -- (4) -- (5) -- (6) -- (7) -- (8) -- (9) -- (10);
		\draw (0) to[out=35,in=145] (4);
		\draw (5) to[out=35,in=145] (8);
		\draw (2) to[out=-35,in=-145] (7);
		\end{scope}
		
		\end{tikzpicture}
		\caption{The graphs $R_{11}$ and $R_{11}'$, both have game domination number $6$.}
		\label{fig:P11}
	\end{center}
\end{figure}
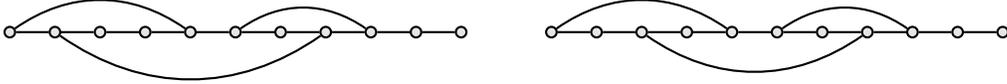

\begin{example}
	\label{exp:path_3edges}
	It holds that $\gamma_g(P_{11}) = 5$, and the same value is achieved for all possibilities of adding one or two edges to the path $P_{11}$. But when three edges are added, it can happen that the game domination number increases to $6 = \lceil \frac{11}{2} \rceil$. Let $V(P_{11}) = [11]_0$ with naturally defined edges. Configurations that result in the game domination number $6$ are obtained by adding the following edges to $P_{11}$: $0\sim4$, $5\sim8$, $1\sim7$ or $0\sim4$, $5\sim8$, $2\sim7$ (see Fig.~\ref{fig:P11}). Denote these graphs with $R_{11}$ and $R_{11}'$, respectively.
\end{example}

By adding the same structure of edges as in the Example~\ref{exp:path_3edges} to the paths $P_{4n+3}$ with vertex set $[4n+3]_0$ for $n \in \{ 2, \ldots, 8\}$, the same phenomena occurs: the value on a path with three additional edges is larger than on the path. It is possible that the same happens for longer paths as well.

Let $n \geq 2$ and let $R_{4n+3}$ be the graph with the vertex set $[4n+3]_0$ obtained from the path $P_{4n+3}$ by adding the edges $0\sim4$, $5\sim8$, and $1\sim7$. Note that $\gamma_g(P_{4n+3}) = 2n + 1$.

\begin{proposition}
	\label{prop:3edges-general}
	If $n \geq 2$, then $\gamma_g(R_{4n+3}) \leq 2n + 2$.
\end{proposition}

\proof
Consider the following subgraphs of $R_{4n+3}$: the graph $R_{11}$ induced by vertices $[11]_0$, and $n-2$ copies of $P_4$ induced by vertices $\{ 4i-1, 4i, 4i+1, 4i+2 \}$ for $i \in \{3, \ldots, n\}$. We have checked by computer that $\gamma_g(R_{11}) = 6$ and that every vertex from $[11]_0$ is an optimal first move for Dominator.

To prove the upper bound, we describe an appropriate strategy of Dominator. His first move is $d_1 = 10$. If Staller plays on $R_{11}$ and this subgraph is not yet dominated, then Dominator selects his optimal move in $R_{11}$. If Staller plays on $R_{11}$ and this subgraph becomes dominated, then Dominator plays his optimal move in one of the subgraphs $P_4$, that is, a vertex of $P_4$ that dominates three vertices of it. If Staller plays on some $P_4$ and it is not yet dominated after her move, then Dominator replies optimally on the same copy of $P_4$ after which all the vertices of the $P_4$ are dominated. Otherwise, Dominator plays his optimal move on some other $P_4$. By the above observation, at most $6$ moves are played on $R_{11}$. As $\gamma_g(P_4) = \gamma_g'(P_4) = 2$, at most $2$ moves are played on each $P_4$. Hence, the number of moves is at most $6 + 2(n-2) = 2n+2$.
\qed

We believe that in Proposition~\ref{prop:3edges-general} the equality actually holds.

Using similar approach, we have also investigated cycles with some additional edges. The following result was obtained by a computer, but no example when the addition of edges would increase the game domination number of a Hamiltonian graph was obtained.
	
\begin{proposition}
	\label{prop:computer_cycles}
	If $4 \leq n \leq 24$ and $G$ is a cycle $C_n$ with two additional edges, then $\gamma_g(G) \leq \gamma_g(C_n) \leq \lceil \frac{n}{2} \rceil$. If $4 \leq n \leq 20$, the same holds for a cycle $C_n$ with three additional edges.
\end{proposition}

\section*{Acknowledgements}

We are grateful to Ga\v{s}per Ko\v{s}mrlj for providing us with his software that computes game domination invariants. We acknowledge the financial support from the Slovenian Research Agency (bilateral grant BI-CN-18-20-008, research core funding P1-0297, projects J1-9109, J1-1693, N1-0095, N1-0108). Kexiang Xu is also supported by NNSF of China (grant No. 11671202) and China-Slovene bilateral grant 12-9.


\newpage
\section*{Appendix just for reviewers}

In Table~\ref{table:double-tadpole} the complete calculations that were used to produce 
Table~\ref{table:exceptions} in the proof of Theorem~\ref{thm:2tailed-tadpole} are listed.

\begin{longtable}{ccc|cc|c}
	$x$ & $y$ & $z$ & $1 + \lceil \w(P'_{4n'+x}) + \w(P'_{4(m'+k')+y+z}) \rceil$ & $\lceil \frac{n(G)}{2} \rceil$ & $\leq$? \\ \hline
	0 & 0 & 0 & $2 k'+2 m'+2 n'+1$ & $2 k'+2 m'+2 n'+2$ & \text{True} \\
	0 & 0 & 1 & $2 k'+2 m'+2 n'+2$ & $2 k'+2 m'+2 n'+2$ & \text{True} \\
	0 & 0 & 2 & $2 k'+2 m'+2 n'+3$ & $2 k'+2 m'+2 n'+3$ & \text{True} \\
	0 & 0 & 3 & $2 k'+2 m'+2 n'+3$ & $2 k'+2 m'+2 n'+3$ & \text{True} \\
	0 & 1 & 0 & $2 k'+2 m'+2 n'+2$ & $2 k'+2 m'+2 n'+2$ & \text{True} \\
	0 & 1 & 1 & $2 k'+2 m'+2 n'+3$ & $2 k'+2 m'+2 n'+3$ & \text{True} \\
	0 & 1 & 2 & $2 k'+2 m'+2 n'+3$ & $2 k'+2 m'+2 n'+3$ & \text{True} \\
	0 & 1 & 3 & $2 k'+2 m'+2 n'+3$ & $2 k'+2 m'+2 n'+4$ & \text{True} \\
	0 & 2 & 0 & $2 k'+2 m'+2 n'+3$ & $2 k'+2 m'+2 n'+3$ & \text{True} \\
	0 & 2 & 1 & $2 k'+2 m'+2 n'+3$ & $2 k'+2 m'+2 n'+3$ & \text{True} \\
	0 & 2 & 2 & $2 k'+2 m'+2 n'+3$ & $2 k'+2 m'+2 n'+4$ & \text{True} \\
	0 & 2 & 3 & $2 k'+2 m'+2 n'+4$ & $2 k'+2 m'+2 n'+4$ & \text{True} \\
	0 & 3 & 0 & $2 k'+2 m'+2 n'+3$ & $2 k'+2 m'+2 n'+3$ & \text{True} \\
	0 & 3 & 1 & $2 k'+2 m'+2 n'+3$ & $2 k'+2 m'+2 n'+4$ & \text{True} \\
	0 & 3 & 2 & $2 k'+2 m'+2 n'+4$ & $2 k'+2 m'+2 n'+4$ & \text{True} \\
	0 & 3 & 3 & $2 k'+2 m'+2 n'+5$ & $2 k'+2 m'+2 n'+5$ & \text{True} \\
	1 & 0 & 0 & $2 k'+2 m'+2 n'+2$ & $2 k'+2 m'+2 n'+2$ & \text{True} \\
	1 & 0 & 1 & $2 k'+2 m'+2 n'+3$ & $2 k'+2 m'+2 n'+3$ & \text{True} \\
	1 & 0 & 2 & $2 k'+2 m'+2 n'+4$ & $2 k'+2 m'+2 n'+3$ & \text{False} \\
	1 & 0 & 3 & $2 k'+2 m'+2 n'+4$ & $2 k'+2 m'+2 n'+4$ & \text{True} \\
	1 & 1 & 0 & $2 k'+2 m'+2 n'+3$ & $2 k'+2 m'+2 n'+3$ & \text{True} \\
	1 & 1 & 1 & $2 k'+2 m'+2 n'+4$ & $2 k'+2 m'+2 n'+3$ & \text{False} \\
	1 & 1 & 2 & $2 k'+2 m'+2 n'+4$ & $2 k'+2 m'+2 n'+4$ & \text{True} \\
	1 & 1 & 3 & $2 k'+2 m'+2 n'+4$ & $2 k'+2 m'+2 n'+4$ & \text{True} \\
	1 & 2 & 0 & $2 k'+2 m'+2 n'+4$ & $2 k'+2 m'+2 n'+3$ & \text{False} \\
	1 & 2 & 1 & $2 k'+2 m'+2 n'+4$ & $2 k'+2 m'+2 n'+4$ & \text{True} \\
	1 & 2 & 2 & $2 k'+2 m'+2 n'+4$ & $2 k'+2 m'+2 n'+4$ & \text{True} \\
	1 & 2 & 3 & $2 k'+2 m'+2 n'+5$ & $2 k'+2 m'+2 n'+5$ & \text{True} \\
	1 & 3 & 0 & $2 k'+2 m'+2 n'+4$ & $2 k'+2 m'+2 n'+4$ & \text{True} \\
	1 & 3 & 1 & $2 k'+2 m'+2 n'+4$ & $2 k'+2 m'+2 n'+4$ & \text{True} \\
	1 & 3 & 2 & $2 k'+2 m'+2 n'+5$ & $2 k'+2 m'+2 n'+5$ & \text{True} \\
	1 & 3 & 3 & $2 k'+2 m'+2 n'+6$ & $2 k'+2 m'+2 n'+5$ & \text{False} \\
	2 & 0 & 0 & $2 k'+2 m'+2 n'+3$ & $2 k'+2 m'+2 n'+3$ & \text{True} \\
	2 & 0 & 1 & $2 k'+2 m'+2 n'+4$ & $2 k'+2 m'+2 n'+3$ & \text{False} \\
	2 & 0 & 2 & $2 k'+2 m'+2 n'+4$ & $2 k'+2 m'+2 n'+4$ & \text{True} \\
	2 & 0 & 3 & $2 k'+2 m'+2 n'+5$ & $2 k'+2 m'+2 n'+4$ & \text{False} \\
	2 & 1 & 0 & $2 k'+2 m'+2 n'+4$ & $2 k'+2 m'+2 n'+3$ & \text{False} \\
	2 & 1 & 1 & $2 k'+2 m'+2 n'+4$ & $2 k'+2 m'+2 n'+4$ & \text{True} \\
	2 & 1 & 2 & $2 k'+2 m'+2 n'+5$ & $2 k'+2 m'+2 n'+4$ & \text{False} \\
	2 & 1 & 3 & $2 k'+2 m'+2 n'+5$ & $2 k'+2 m'+2 n'+5$ & \text{True} \\
	2 & 2 & 0 & $2 k'+2 m'+2 n'+4$ & $2 k'+2 m'+2 n'+4$ & \text{True} \\
	2 & 2 & 1 & $2 k'+2 m'+2 n'+5$ & $2 k'+2 m'+2 n'+4$ & \text{False} \\
	2 & 2 & 2 & $2 k'+2 m'+2 n'+5$ & $2 k'+2 m'+2 n'+5$ & \text{True} \\
	2 & 2 & 3 & $2 k'+2 m'+2 n'+6$ & $2 k'+2 m'+2 n'+5$ & \text{False} \\
	2 & 3 & 0 & $2 k'+2 m'+2 n'+5$ & $2 k'+2 m'+2 n'+4$ & \text{False} \\
	2 & 3 & 1 & $2 k'+2 m'+2 n'+5$ & $2 k'+2 m'+2 n'+5$ & \text{True} \\
	2 & 3 & 2 & $2 k'+2 m'+2 n'+6$ & $2 k'+2 m'+2 n'+5$ & \text{False} \\
	2 & 3 & 3 & $2 k'+2 m'+2 n'+6$ & $2 k'+2 m'+2 n'+6$ & \text{True} \\
	3 & 0 & 0 & $2 k'+2 m'+2 n'+3$ & $2 k'+2 m'+2 n'+3$ & \text{True} \\
	3 & 0 & 1 & $2 k'+2 m'+2 n'+4$ & $2 k'+2 m'+2 n'+4$ & \text{True} \\
	3 & 0 & 2 & $2 k'+2 m'+2 n'+5$ & $2 k'+2 m'+2 n'+4$ & \text{False} \\
	3 & 0 & 3 & $2 k'+2 m'+2 n'+5$ & $2 k'+2 m'+2 n'+5$ & \text{True} \\
	3 & 1 & 0 & $2 k'+2 m'+2 n'+4$ & $2 k'+2 m'+2 n'+4$ & \text{True} \\
	3 & 1 & 1 & $2 k'+2 m'+2 n'+5$ & $2 k'+2 m'+2 n'+4$ & \text{False} \\
	3 & 1 & 2 & $2 k'+2 m'+2 n'+5$ & $2 k'+2 m'+2 n'+5$ & \text{True} \\
	3 & 1 & 3 & $2 k'+2 m'+2 n'+5$ & $2 k'+2 m'+2 n'+5$ & \text{True} \\
	3 & 2 & 0 & $2 k'+2 m'+2 n'+5$ & $2 k'+2 m'+2 n'+4$ & \text{False} \\
	3 & 2 & 1 & $2 k'+2 m'+2 n'+5$ & $2 k'+2 m'+2 n'+5$ & \text{True} \\
	3 & 2 & 2 & $2 k'+2 m'+2 n'+5$ & $2 k'+2 m'+2 n'+5$ & \text{True} \\
	3 & 2 & 3 & $2 k'+2 m'+2 n'+6$ & $2 k'+2 m'+2 n'+6$ & \text{True} \\
	3 & 3 & 0 & $2 k'+2 m'+2 n'+5$ & $2 k'+2 m'+2 n'+5$ & \text{True} \\
	3 & 3 & 1 & $2 k'+2 m'+2 n'+5$ & $2 k'+2 m'+2 n'+5$ & \text{True} \\
	3 & 3 & 2 & $2 k'+2 m'+2 n'+6$ & $2 k'+2 m'+2 n'+6$ & \text{True} \\
	3 & 3 & 3 & $2 k'+2 m'+2 n'+7$ & $2 k'+2 m'+2 n'+6$ & \text{False} \\
	
	\caption{For all possible values of $x, y, z$ we check if\\ $1 + \lceil \w(P'_{4n'+x}) + \w(P'_{4(m'+k')+y+z}) \rceil \leq \lceil \frac{n(G)}{2} \rceil$.}
	\label{table:double-tadpole}
\end{longtable}

\end{document}